\documentclass[preprint,10pt,authoryear]{elsarticle}

\usepackage[a4paper, top=2.5cm, bottom=2.5cm, left=2.5cm, right=2.5cm]{geometry}
\usepackage{setspace}
\usepackage{amssymb}
\usepackage{amsmath}
\usepackage{comment}
\usepackage{caption,subcaption}     
\usepackage[figuresright]{rotating} 
\usepackage{algorithmic}            
\usepackage{algorithm}              

\usepackage{multirow}               
\usepackage[table,xcdraw]{xcolor}   
\usepackage{arydshln}               

\usepackage{amsthm}                 

\journal{ArXiv}

\begin{document}

\begin{frontmatter}

\title{Mixed-model Sequencing with Reinsertion of Failed Vehicles: A Case Study for Automobile Industry}
\tnotetext[t1]{This research did not receive any specific grant from funding agencies in the public, commercial, or
not-for-profit sectors.}

\author[1]{I. Ozan Yilmazlar\corref{cor1}}
\ead{iyilmaz@clemson.edu}

\author[1]{Mary E. Kurz}
\ead{mkurz@clemson.edu}

\cortext[cor1]{Corresponding author}
\address[1]{Department of Industrial Engineering, Clemson University, Clemson SC, 29634}

\begin{abstract}
    In the automotive industry, some vehicles, failed vehicles, cannot be produced according to the planned schedule due to some reasons such as material shortage, paint failure, etc. These vehicles are pulled out of the sequence, potentially resulting in an increased work overload. On the other hand, the reinsertion of failed vehicles is executed dynamically as suitable positions occur. In case such positions do not occur enough, either the vehicles waiting for reinsertion accumulate or reinsertions are made to worse positions by sacrificing production efficiency.
    
    This study proposes a bi-objective two-stage stochastic program and formulation improvements for a mixed-model sequencing problem with stochastic product failures and integrated reinsertion process. Moreover, an evolutionary optimization algorithm, a two-stage local search algorithm, and a hybrid approach are developed. Numerical experiments over a case study show that while the hybrid algorithm better explores the Pareto front representation, the local search algorithm provides more reliable solutions regarding work overload objective. Finally, the results of the dynamic reinsertion simulations show that we can decrease the work overload by ~20\% while significantly decreasing the waiting time of the failed vehicles by considering vehicle failures and integrating the reinsertion process into the mixed-model sequencing problem.
\end{abstract}

\end{frontmatter}

\section{Introduction}\label{section: reinsertion-introduction}
An assembly line is a manufacturing system in which a product is progressively built by moving it through a sequence of workstations, each specializing in a specific task.
Mixed-model assembly lines (MMAL) are developed to produce multiple product variations or models on the same assembly line, allowing for flexibility and customization while maximizing efficiency and productivity.
The variety of models produced significantly increases as the product becomes more complex and customizable. An automobile is a useful example of such a product, e.g., a German car manufacturer has up to $10^{24}$ potential configurations for their vehicles \cite{pil2004linking}. Some products may require additional steps, specialized equipment, or longer assembly processes due to their unique attributes, resulting in substantial variations in processing times. The tasks are distributed through the assembly line balancing, ensuring that the average workload of each station aligns with the established \textit{cycle time}, the fixed duration between two consecutive product launches. However, a poor sequence of products may result in uneven workloads across workstations, excessive work allocation (\textit{work overload}), or idle time, due to the varying processing times of products. When such a work overload situation occurs, interventions are needed such as stopping the line or utilizing more workers. In order to avoid or minimize such interventions, the assembly process is balanced by sequencing the products. The corresponding sequencing problem is called as a mixed-model sequencing (MMS) problem which minimizes the work overload duration across all workstations.

In the automobile industry, car manufacturers have adjusted their assembly lines to accommodate the mixed-model production of vehicles with diesel and gasoline engines, however, the starting to produce electric vehicles (EV) (or hybrid vehicles) on the same line has posed new challenges. In contrast to other vehicles, EVs have large batteries, resulting in substantial differences in tasks, particularly at the station where the battery is loaded. Accordingly, the position of EVs in a sequence needs to be carefully determined since too many EVs in a subsequence may result in a significant work overload at the corresponding stations. 
To illustrate further, consider an MMAL where a station is responsible for an electric battery system. The processing time for installing a battery system on an EV might take, on average, three minutes, while tasks on a non-EV might take around 30 seconds. Without proper sequencing, the battery loading station could experience work overloads and idle times due to the drastic processing time difference. By solving the MMS problem for each planning horizon, the workload can be distributed in a way that accommodates the varying processing times, ensuring efficient utilization of resources and maximizing production output.

As car manufacturers increase the proportion of the EVs in the product mix in response to market demand, the issues that cause a vehicle not to be produced in the planned sequence become more significant. For instance, if a part is delayed from a supplier, or if there is a paint issue with a vehicle, a planned sequence of vehicles may have a missing vehicle, which is referred to as \textit{failed vehicles}. Even if this missing vehicle has a gasoline or diesel engine, its absence can still affect the stations that require intense battery-related tasks. The failed vehicles are pulled out of the sequence and the resulting gap is closed either by moving the succeeding vehicles forward or by filling with another vehicle that failed previously and waiting to be reinserted, we call these vehicles, failed but ready for reinsertion, as \textit{reinstating vehicles} throughout this study. In a traditional MMS, the failure of vehicles is not considered which means that the failure of vehicles may result in a significant additional work overload. 

Moreover, in automobile manufacturing facilities, the reinsertion process is executed dynamically. Ideally, a failed vehicle that is ready for reinsertion is reinserted into the sequence only if an appropriate position is found, i.e., reinsertion is executed only if it does not cause an extra work overload. The vehicles to be reinserted can accumulate if appropriate positions do not occur enough. Once the number of reinstating vehicles exceeds a threshold, the line is stopped and all the reinstating vehicles are produced back to back which generally impacts production efficiency. As a result, we propose integrating the possibility of reinsertion into the sequencing model in order to minimize the possible work overload due to vehicle failures and reinsertions and to avoid line stoppages due to the excessive number of reinstating vehicles. 

\begin{figure}[h]
\captionsetup{justification=centering}
\centering
\begin{subfigure}{0.5\textwidth}
  \centering
  \includegraphics[width=0.85\linewidth]{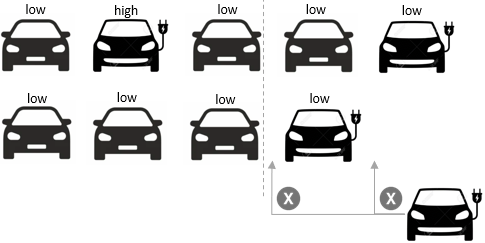}
  \caption{non-robust sequence}
  \label{fig_sub: reinsertion-non-robust sequenc}
\end{subfigure}%
\begin{subfigure}{0.5\textwidth}
  \centering
  \includegraphics[width=0.85\linewidth]{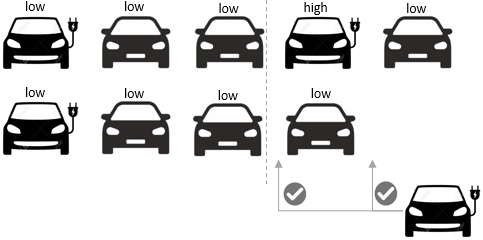}
  \caption{robust sequence}
  \label{fig_sub: reinsertion-robust sequence}
\end{subfigure}
\caption{Illustration of a non-robust and robust sequence to stochastic failures and reinsertion process}
\label{fig: reinsertion-robust-sequence}
\end{figure}

In Figure \ref{fig: reinsertion-robust-sequence}, we motivate a robust sequence that considers vehicle failures and the reinsertion process of failed vehicles for a battery installation station where back-to-back EVs cause a large amount of work overload. In the example, there are five vehicles to be sequenced, three non-EVs and two EVs, and one of the EVs has a high failure probability. There is one failed vehicle that will be ready for reinsertion starting from position 3. Assuming the high-risk vehicle fails, different consequences occur based on the initial sequence. In the non-robust sequence, the reinsertion of the reinstating vehicle results in work overload for either of the possible positions. On the other hand, in the robust sequence, the reinsertion to any possible position does not cause any work overload. 

This study aims to generate robust sequences that take into account the potential vehicle failures and their reinsertion back into the sequence, in order to avoid additional workloads. We focus on the final assembly line, where we assume that vehicles must follow the original sequence they arrive in from the paint shop, without the option of resequencing at the point where the vehicle failures are realized. Accordingly, when a vehicle is taken out of the sequence, the subsequent vehicles fill the resulting gap. The key contributions of this study are as follows:

\begin{itemize}
\item A bi-objective two-stage stochastic program is presented for an MMS problem with stochastic product failures and reinsertion process, and formulation improvements are provided. This is the first study that considers the product reinsertion process in MMS, to the best of our knowledge.
\item We develop a bi-objective evolutionary algorithm, a two-stage bi-objective local search algorithm, and a local search integrated bi-objective evolutionary algorithm to tackle the proposed problem. The numerical results show that the local search integrated evolutionary algorithm outperforms other algorithms by providing work reliable solutions in terms of work overload objective. On the other hand, the hybrid local algorithm provides a better exploration of solution space
\item A case study with the data inspired by our car manufacturer industrial partner is conducted via dynamic reinsertion simulations. The numerical experiments show that the work overload can be reduced by around 20\% by considering stochastic car failures and the reinsertion process, while significantly decreasing the waiting time of failed vehicles until their reinsertion.
\end{itemize}

The remainder of this paper is structured as follows. The related literature is reviewed in Section \ref{section: reinsertion-literature}. The proposed problem is defined, and mathematically formulated in Section \ref{section: reinsertion-problem-definition}. Efficient solution approaches to tackle the problem are presented in Section \ref{section: reinsertion-solution-approach}. The numerical experiments that analyze the performance of presented solution approaches are executed, and the results are shared in Section \ref{section: reinsertion-numerical-experiments}. Finally, a summary of our findings and a discussion about future work are presented in Section \ref{section: reinsertion-conclusion}.

\section{Related Literature}\label{section: reinsertion-literature}
The efficiency of an MMAL primarily relies on the management of two key decision problems \citep{bard1992analytic}. As a long-term decision problem, line balancing has a strategic role, as it involves allocating the workload across the workstations in the assembly line \citep{becker2006survey, boysen2007classification, sivasankaran2014literature, pearce2019effective, yilmazlar2020case}. As a short-term decision problem, the sequencing problem has an operational role, as it addresses daily production planning. Combining mixed-model assembly line balancing and sequencing problems has garnered significant interest among researchers since considering strategic and operational decisions together provides benefits to the line efficiency \citep{ozcan2010balancing, lopes2020iterative, peng2022research}. For the sequencing problem, researchers have developed two main categories of objective functions: those related to work overload and those focused on just-in-time objectives. Additionally, the literature presents three solution approaches: level-scheduling \citep{steiner1993level}, car-sequencing \citep{dincbas1988solving, solnon2008car}, and MMS. A comprehensive survey, including classification schemes and detailed explanation, is given by \citet{boysen2009sequencing}.

The problem proposed in this paper can be considered a combined mixed-model sequencing and resequencing problem. The resequencing problem minimizes the work overload while reinserting a set of failed vehicles into a given sequence \citep{boysen2012resequencing}. There are several studies in the literature that tackle the resequencing problem in the automobile industry. \citet{franz2014resequencing} considers a resequencing setting in the automotive industry. They compare the performance of three local search heuristic features approaches while minimizing the total work overload occurrences. A variable neighborhood search combined with tabu search outperforms other feature combinations. \citet{franz2015dynamic} considers the resequencing problem in a dynamic setting, dynamically supplied failed vehicles, and compares two reinsertion strategies integrated into local search heuristics; reinsert the vehicle as soon as it is ready and wait for a good reinsertion position. They show that both strategies perform better, compared to real-world application that does not utilize such strategy. \citet{gujjula2009resequencing} considers the resequencing problem in the level-scheduling setting. \citet{gunay2017stochastic} proposes a two-stage stochastic model for a resequencing problem that considers the post-paint restoring buffer, which is used to restore the initial sequence before the paint quality issues. \citet{leng2022multi} proposes a multi-objective reinforcement learning approach to solve combined resequencing and color-batching problems.

While most of the existing literature on MMS primarily focuses on models that have deterministic parameters, there are a limited number of studies that examine stochastic parameters, only related to processing times or demand. \citet{zhao2007modeling} introduces a Markov chain-based approach that minimizes the expected total work overload duration, by generating sub-intervals that represent possible positions of workers within the stations based on the stochastic processing times. \citet{mosadegh2017heuristic} proposes a heuristic approach inspired by Dijkstra's algorithm to address a single-station MMS with stochastic processing times, formulated as a shortest path problem. \citet{mosadegh2020stochastic} formulates a multiple-station MMS with stochastic processing times as a mixed-integer linear programming (MILP) problem. They introduce a Q-learning-based simulated annealing heuristic to solve real-world-sized problems and demonstrate a decrease in expected work overload compared to the deterministic problem. \citet{brammer2022stochastic} suggests a reinforcement learning approach for solving MMS with stochastic processing times. Their results indicate that the proposed method outperforms SA and GA by at least 7\% in terms of solution quality. Furthermore, stochastic parameters are also considered in integrated mixed-model balancing and sequencing problems. Specifically, stochastic processing times are studied in \citep{agrawal2008collaborative, ozcan2011genetic, dong2014balancing} and stochastic demand is studied in \citet{sikora2021benders}.

To the best of our knowledge, there is currently no existing research that focuses on establishing robust sequences that handle work overloads resulting from product failures by considering the reinsertion process within any sequencing structure. The only studies that consider the stochastic product failure, but not the reinsertion of failed products, are \citet{hottenrott2021robust} within car sequencing structure, and \citet{yilmazlar2023mixed} within MMS structure.

\section{Problem Statement and Mathematical Formulation}\label{section: reinsertion-problem-definition}
In Section \ref{section: reinsertion-problem-statement}, we define the MMS with reinsertion and illustrate the problem with an example. Then, in Section \ref{section: reinsertion-mathematical-model}, we provide a mixed-integer quadratic program for the proposed problem.

\subsection{Problem Statement}\label{section: reinsertion-problem-statement}
In an MMAL (Mixed-Model Assembly Line), a conveyor belt is used to interconnect a group of workstations. Products are launched at a fixed rate, and the duration between two consecutive products is referred to as the cycle time $c$. The products move along the belt with a constant speed of one time unit (TU). We assume that the length of station $k\in K$ is $l_{k} \geq c$ in TU, which represents the total time required for a workpiece to traverse the station. 
Operators are responsible for carrying out assigned tasks within the station length, i.e., closed-border stations. If they fail to complete their tasks within the allocated time, a utility worker intervenes to complete the remaining work. The additional workload that remains unfinished is known as \textit{work overload}. The efficiency of the assembly line is greatly impacted by the sequence of products. Therefore, MMS  determines the sequence of a given set of products $V$ by assigning each product $v\in V$ to a specific position $t\in T$.

In this study, we consider that each vehicle is unique since our automobile manufacturer partner produces vehicles based on custom orders and customization offers billions of different options. Hence, we sequence vehicles instead of configurations. Accordingly, we define the first-stage binary decision variable $x_{vt}$, which takes the value of  1 if vehicle $v \in V$ is assigned to position $t \in T$. Additionally, the second-stage decision variables are defined to determine the reinsertion of failed vehicles and the final position of all vehicles after the fails are realized and the reinsertions are decided. The binary decision variable $y_{it}$ is 1 if failed vehicle $i\in F$ reinserted at position $t$, and binary decision variable $s_{vt}$ is 1 if the final position of vehicle $v\in V\cup F$ is $t$. 
We denote the processing time of vehicle $v \in V$ at station $k \in K$ by $p_{kv}$. The starting position and work overload of the vehicle at position $t=1,\ldots,|T|+|F_{old}|$ for station $k \in K$ are denoted by $z_{kt}$ and $w_{kt}$, respectively. In Table \ref{tbl: reinsertion-notations}, we present all the sets, parameters, and decision variables used in the proposed mathematical formulation. The second-stage decision variables are scenario-dependent, however, such dependency is dropped for the sake of notation simplicity.

\begin{table}[]
\captionsetup{format=plain, labelfont=bf, justification=centering, singlelinecheck=false, font=small, skip=4pt}
\caption{List of parameters and decision variables used in the model}
\centering
\scalebox{0.75}{
\begin{tabular}{
>{\columncolor[HTML]{FFFFFF}}l 
>{\columncolor[HTML]{FFFFFF}}l lll}
\textbf{Sets and Index}     & { \textbf{}}                                                                               &  &  &  \\ \cline{1-2}
$V, v$                          & Vehicles that are initially planned to be produced in the current horizon               &  &  &  \\
$F_{new}, i$                    & Failed vehicles that were planned to be produced in the current horizon          &  &  &  \\
$F_{old}, i$                    & Failed vehicles that were planned to be produced in a previous horizon       &  &  &  \\
$K, k$                          & Stations                                                                                      &  &  &  \\
$T, t$                          & Positions                                                                                     &  &  &  \\
$\Omega, \omega$                          & Scenarios                                                                           &  &  &  \\
                              &                                                                                               &  &  &  \\
\textbf{Parameters}           & \textbf{}                                                                                     &  &  &  \\ \cline{1-2}
$p_{kv}$                        & The processing time of vehicle $v$ at station $k$                                           &  &  &  \\
$l_k$                           & The length of station $k$                                                                        &  &  &  \\
$c$                             & The cycle time                                                                                &  &  &  \\
$f_v$                           & The failure probability of vehicle $v$                                                          &  &  &  \\
$e_{vn}$                        & 1 if vehicle $v$ exists at scenario $n$, 0 otherwise                                          &  &  &  \\
$f_{max}$                       & The maximum number of failed cars allowed at the end of the planning horizon                     &  &  &   \\
$g_{i}$                         & The number of days between the planned and current period for the failed car $i\in F_{old}$      &  &  &   \\
$d_{i}$                         & The number of extra days until the delivery date that the failed car $i\in F_{old}$ had when failed &  &  &  \\
$\lambda$                       & The minimum number of positions between each reinsertion                                   &  &  &   \\
$r_{i}$                         & The number of positions required for failed vehicle $i \in F$ to be ready for reinsertion                     &  &  &   \\
                              &                                                                                               &  &  &  \\
\textbf{First-Stage Decision Variables} & \textbf{}                                                     &  &  &  \\ \cline{1-2}
$x_{vt}$                           & 1 if vehicle $v \in V$ is assigned to position $t \in T$, 0 otherwise                &  &  &  \\
$\delta_{vt}$       & 1 if vehicle $v\in V$ is positioned after position $t$ before reinsertions, 0 otherwise      &  &  &  \\
$\delta_{vt}^{'}$   & 1 if vehicle $v\in V$ is positioned at or after position $t$ before reinsertions, 0 otherwise      &  &  &  \\
                      &  &  &  \\
\textbf{Second-Stage Decision Variables} & \textbf{}                                                           &  &  &  \\ \cline{1-2}
$w_{kt}$            & The work overload at station $k \in K$ at cycle $t \in \Bar{T}$                                &  &  &  \\
$z_{kt}$            & Starting position of operator at station $k \in K$ at the beginning of cycle $t \in \Bar{T}$   &  &  &  \\
$b_{kt}$            & The processing time at station $k \in K$ at cycle $t \in \Bar{T}$                              &  &  & \\
$y_{it}$            & 1 if failed vehicle $i\in F$ is inserted to position $t=\{1,\ldots,|T|+1\}$, 0 otherwise                    &  &  &  \\
$s_{vt}$            & 1 if vehicle $i\in \Bar{V}$ is sequenced at position $t \in \Bar{T}$ after reinsertions, 0 otherwise &  &  &  \\
$\gamma_{it}$       & 1 if failed vehicle $i\in F$ is inserted after position $t \in T$, 0 otherwise      &  &  & 
\end{tabular}}
\label{tbl: reinsertion-notations}
\end{table}

There are two types of failed vehicles that require special consideration for the reinsertion process. The first type is the vehicles that had failed in a previous production horizon $F_{old}$, have not been reinserted yet, and become ready for reinsertion at the beginning or during the current horizon. The second type is the vehicles that fail in the current production horizon $F_{new}$. In addition to the sets defined in Table \ref{tbl: reinsertion-notations}, we define some other useful sets, based on the relation of the already defined sets. The set $F=F_{new}\cup F_{old}$ presents the set of all failed vehicles. The set $V^{'}=V\setminus F_{new}$ denotes the vehicles planned for the current production horizon and do not fail. The sets $\Bar{V}=V\cup F_{old}$ and $\Bar{T}$ denote the set of all vehicles and  positions after the reinsertions, respectively. 

In this study, while we adopt the basic assumptions of the MMS problem as given by Bolat \emph{et al.} \cite{bolat1992scheduling}, additional rules/assumptions for the vehicle failures and the reinsertion process are defined as follows:
\begin{itemize}
    \item  It is assumed that each vehicle $v \in V$ has a failure probability $f_v$, and failures are independent of each other. In our numerical experiments in Section \ref{section: reinsertion-numerical-experiments}, the failure probabilities are estimated from the historical data by doing feature analysis and using logistic regression.
    \item The vehicles go through the body and paint shops in the scheduled sequence before the assembly process. Hence, a failed vehicle must be pulled out of the sequence, however, its position can be compensated by a reinstating vehicle.
    \item Once the failure of a vehicle occurs due to any reason, there is an uncertain amount of time (measured in the number of positions) required for the failure reason to be fixed, i.e., the number of positions required for a failed vehicle to become a reinstating vehicle. For example, repainting of the vehicle if the issue was paint quality, or arrival of material if it was a supply chain issue. This time, referred to as \textit{being ready time}, is modeled as a discrete uniform distribution in a closed range. 
    \item There is an upper limit on the number of reinstating vehicles at the end of each production horizon, $f_{max}$. In the car manufacturing facilities, if $f_{max}$ is to be exceeded, the whole line stops and only produces the reinstating vehicles. The motivation behind this application is to avoid the late delivery of failed vehicles.
    \item The number of failed vehicles at the beginning of a production horizon is predicted randomly using a discrete uniform distribution between 0 and $f_{max}$.
    \item It is assumed that only one reinsertion can be done into each subsequence of length $\lambda$.
    \item A failed vehicle from a previous planning horizon must be reinserted if the current horizon is the last day for it to be started, $g_{i} = d_{i}$, so that the vehicle can be delivered to the customer on time.
\end{itemize}

The integration of the reinsertion process into the sequencing problem introduces conflicting objectives. The number of vehicles reinserted during a planning horizon conflicts with the total work overload since the reinsertions are made into an already optimized sequence, minimized work overload. Throughout the paper, these objectives are referred to as \textit{work overload objective} and \textit{reinsertion objective}, respectively. Accordingly, we employ two objectives for our problem: minimizing total work overload duration and minimizing the sum of squared waiting days of each not reinserted failed vehicle. 
For the work overload handling procedure, we adopt the side-by-side policy which assumes that the regular operator stops working on the workpiece once the workpiece reaches the station border. The remaining job of the workpiece is completed by the so-called utility worker, and the regular operator starts working on the next workpiece at position $l_k - c$ in the same station. 

\subsection{Mathematical Model Under Uncertainty}\label{section: reinsertion-mathematical-model}
The operational dynamics of the reinsertion process, which is executed after the realization of vehicle failures, motivated us to  formulate our problem as a two-stage stochastic program. In the first stage (here-and-now), the vehicle sequence is decided before any car failures are realized. Subsequently, when the car failures become known, the work overload is minimized through the second-stage decisions (wait-and-see) based on the initial sequence. While the determination of first-stage decisions involves assigning each vehicle to a specific position with the objective of minimizing the anticipated objective function value in the second stage, the second-stage decisions involve the reinsertions with the work overload and reinsertion objectives. For the sake of clarity, throughout the paper, we call the decisions on either a failed vehicle reinserted or not as \textit{binary reinsertion decisions}, on the other hand, the decision on the reinsertion positions as \textit{reinsertion position decisions}.

In the proposed problem, the uncertainty lies in several factors: 1. vehicle failures in the current horizon, each vehicle either exists or fails, 2. the number of vehicles in $F_{old}$, 3. the number of days between the planned and current period $g_{i}$ for each failed car $i\in F_{old}$, 4. The number of positions required for each failed vehicle $i\in F$ to be ready for reinsertion $r_{i}$. There are a total of $2^{|V|}$, $(f_{max}+1)$, $G^{|F_{old}|}$ (assume that $G$ is the lead time of the facility), and $|V|*|F|$ scenarios for the factors 1, 2, 3, and 4, respectively. In order to decrease the degree of scenario space, we predict factors 3 and 4 and fix them as deterministic parameters. That is, $g_{i}$ and $r_{i}$ are known beforehand and used in case of vehicle $i$ failure. Accordingly, we have a total of $2^{|V|}*(f_{max}+1)$ scenarios. 

To formulate the problem, we represent varying realizations of vehicle failures and failed vehicles by a set of finite scenarios $\Omega$.  Each scenario $\omega \in \Omega$ is denoted by the vehicle failure, let $e_{v\omega}=0$ if vehicle $v$ fails and $e_{v\omega}=1$ if vehicle $v\in V$ exists at scenario $\omega \in \Omega$, and the number of failed vehicles from previous horizons $|F_{old}|$ is uniformly distributed in $[0, f_{max}]$. The probability of scenario $\omega$ can be calculated as $p_{w}=\frac{1}{f_{max}+1}\prod_{v=1}^{|V|}f_v^{1-e_{v\omega}} (1-f_v)^{e_{v \omega}}$.

We note that the sequence length depends on the vehicle failures and the second-stage decisions (reinsertions), and such a problem cannot be solved using state-of-art commercial solvers. Hence, we develop some methods to improve the formulation so that we can formulate the second-stage as a mixed-integer quadratically constrained program (MIQCP). 
First, we create a \textit{dummy position} at $|T|+1$ to where we make dummy reinsertions, in other words, we assume that all failed vehicles $i\in F$ that are actually not reinserted are assigned to position $|T|+1$. This trick helps each scenario have a fixed sequence length that is no longer decision dependent. The sequence length of scenario $\omega$ is $|T|+|F_{old}^{\omega}|$.
Next, dummy reinsertions, that are made at the end of the sequence, should be exempt from the work overload calculation since they actually do not exist in the final sequence (not reinserted), and they will be produced in a future horizon. In order to tackle this issue, we turn each vehicle that is reinserted at the dummy position into a \textit{neutral vehicle} by setting the processing time equal to the cycle for each station. We refer to these vehicles as neutral vehicles because they do not have any impact on the schedule in terms of work overload. A n extensive explanation of neutral vehicles is given by \citet{yilmazlar2023mixed}.  

Finally, we present a two-stage stochastic program for the \textit{full-information problem} where all possible realizations are considered.

\begin{subequations}
\label{formulation: reinsertion-first-stage}
\begin{flalign}
\min_{x,\delta,\delta^{'}} \quad & \sum_{\omega\in \Omega}\rho_\omega Q((x,\delta,\delta^{'}),\omega) \label{c1.0.1} \\
\text{s.t.} \quad & \sum_{v\in V}x_{vt} = 1, \qquad \  t\in T & \label{c1.1.1} \\ 
& \sum_{t\in T}x_{vt} = 1, \qquad \  v\in V & \label{c1.1.2} \\ 
& \sum_{t\in T}\delta_{vt} - x_{vt}t = 0, \qquad \  v\in V & \label{c1.6.1} \\
& \delta_{vt} - \delta_{v(t-1)} \leq 0, \qquad \  v\in V, \enspace \  t=2,\ldots,|T| & \label{c1.6.2} \\
& \delta_{v0}^{'} = 0, \qquad v\in V & \label{c1.7.1} \\
& \delta_{vt}^{'} - \delta_{v(t-1)} = 0, \qquad v\in V, \enspace t=2,\ldots, |T| & \label{c1.7.2} \\
& x \in \{0,1\}, \label{c1.3.1} \\
& \delta, \delta^{'} \geq 0 & \label{c1.3.2} 
\end{flalign}
\end{subequations}
where $Q(x,\omega) =$
\begin{subequations}
\label{formulation: reinsertion-second-stage}
\begin{flalign}
\min_{y,s,\gamma,z,w,b,\beta} \quad & \sum_{t\in \Bar{T}}\sum_{k\in K} w_{kt} + \sum_{i\in F}(1-\sum_{t\in T}y_{it})(g_{i}+1)^2  \label{c1.0.2} \\
\text{s.t.} \quad & x_{it} - \sum_{h=\min(|T|, t+r_{i})}^{|T|+1}y_{ih} \leq 0 \qquad i\in F_{new}, \enspace t\in T & \label{c1.2.1} \\ 
& \sum_{t=1}^{|T|+1}y_{it}t \geq r_{i}, \qquad i\in F_{old} & \label{c1.2.2} \\ 
& \sum_{t\in T}y_{it} = 1, \qquad i\in \{F_{old} | g_{i} = d_{i}\} & \label{c1.2.3} \\
& \sum_{h=t}^{t+\lambda-1}\sum_{i\in F}y_{ih} \leq 1, \qquad  t= 1,\ldots,|T|-\lambda+1 & \label{c1.2.4} \\
& \sum_{t\in T}(1-y_{it}) \leq f_{max}, \qquad i\in F & \label{c1.2.5} \\
& \sum_{t=1}^{|T|+1}y_{it} = 1, \qquad i\in F & \label{c1.5} \\
& \sum_{t\in T}\gamma_{it} - y_{it}t = 0, \qquad \  i\in F  & \label{c1.8.1} \\
& \gamma_{it} - \gamma_{i(t-1)} \leq 0, \qquad \  i\in F, \enspace \  t=2,\ldots,|T| & \label{c1.8.2} \\
& \sum_{v\in \Bar{V}}s_{vt} = 1, \qquad \  t\in \Bar{T} & \label{c1.9.1} \\ 
& \sum_{t\in \Bar{T}}s_{vt} = 1, \qquad \  v\in \Bar{V} & \label{c1.9.2} \\ 
& \sum_{t\in \Bar{T}} s_{vt}t - \sum_{t\in T}x_{vt}t - \sum_{t\in T}\sum_{i\in F}y_{it}\delta_{vt}^{'} + 
\sum_{i\in F_{new}} \sum_{t\in T}x_{it}\delta_{vt}(1-e_{i}) = 0, \qquad \  v\in V^{'} & \label{c1.10.1} \\ 
& \sum_{t\in \Bar{T}} s_{it}t - \sum_{t\in T}y_{it}t - \sum_{t\in T}\sum_{j\in F}y_{jt}\gamma_{it} + 
\sum_{j\in F_{new}} \sum_{t\in T}x_{jt}\gamma_{it}(1-e_{j}) - y_{vT}\beta_{i} = 0,  \qquad \  i\in F & \label{c1.10.2} \\ 
& b_{kt} - \sum_{v\in V^{'}} p_{kv}s_{vt} - \sum_{i\in F}s_{it}(p_{ki}(1-y_{i(|T|+1)}) + cy_{i(|T|+1)}) = 0, \qquad \  k\in K, \enspace \  t\in \Bar{T} & \label{c1.20} \\
&  z_{kt}+b_{kt}-z_{k(t+1)}-w_{kt} \leq c, \qquad \  k\in K, \enspace \   t\in \Bar{T} & \label{c1.21} \\ 
&  z_{kt}+b_{kt}-w_{kt} \leq l_k, \qquad \  k\in K, \enspace \   t\in \Bar{T} & \label{c1.22} \\ 
&  z_{k0} = 0, \qquad \  k\in K & \label{c1.23} \\ 
&  z_{k(|\Bar{T}|+1)} = 0, \qquad \  k\in K & \label{c1.24} \\ 
& y, s, \gamma \in \{0,1\},  & \label{c1.25} \\
& z, w, b, \beta \geq 0 & \label{c1.26} 
\end{flalign}
\end{subequations}


The first-stage problem \eqref{formulation: reinsertion-first-stage} minimizes the expected cost associated with the second-stage problem.
Constraint sets \eqref{c1.1.1} and \eqref{c1.1.2} ensure that each position has one vehicle assigned and each vehicle is assigned to one position, respectively. 
Constraint sets \eqref{c1.6.1} and \eqref{c1.6.2} builds the relationship between $x$ and $\delta$ variables. The $\delta$ variables help to determine the number of failed vehicles up to each position in constraint set \eqref{c1.10.1}.
Constraint sets \eqref{c1.7.1} and \eqref{c1.7.2} build the relationship between $\delta^{'}$ and $\delta$ variables. The $\delta^{'}$ variables help to determine the number of reinserted vehicles up to each position in constraint set \eqref{c1.10.1}.
Constraint sets \eqref{c1.3.1} and \eqref{c1.3.2} represent the domain of the first-stage variables.

In the second-stage problem \eqref{formulation: reinsertion-second-stage}, the objective function minimizes the sum of the total work overload duration and the sum of the square of the number of waiting days of the failed vehicles that arenot to be reinserted, given the first-stage decision (sequence) and scenario $\omega \in \Omega$.  
Constraint sets \eqref{c1.2.1} - \eqref{c1.2.5} force the restrictions about the reinsertion. The constraint sets \eqref{c1.2.1} and \eqref{c1.2.2} ensure that new and previously failed vehicles are not reinserted before they are ready, respectively. Constraint set \eqref{c1.2.3} guarantees that a failed vehicle at its due must be reinserted. The constraint set \eqref{c1.2.4} assures that each subsequence in length $\lambda$ has at most one reinsertion. Constraint set \eqref{c1.2.5} assures that the limit on the maximum number of not reinserted vehicles is not exceeded.
Constraint set \eqref{c1.5} ensures that all failed vehicles are reinserted; the reinsertions at position $|T|+1$ are the dummy reinsertions. 
Constraint sets \eqref{c1.8.1} and \eqref{c1.8.2} build the relationship between $x$ and $\gamma$ variables. The $\gamma$ variables help to determine the number of failed and reinserted vehicles up to each position in constraint set \eqref{c1.10.2}.
Constraint sets \eqref{c1.9.1} and \eqref{c1.9.2} ensure that, after the reinsertions, each position has one vehicle assigned and each vehicle is assigned to one position, respectively. Note that $\Bar{V}$ and $\Bar{T}$ include the previously failed vehicles and additional positions for them, respectively. 
Constraint sets \eqref{c1.10.1} and  \eqref{c1.10.2} determine the final position of each non-failed and failed vehicle, respectively. Note that $V^{'}$ denotes the vehicles planned for the current production horizon and do not fail. In constraint set \eqref{c1.10.1}, the final position of a non-failed vehicle is set to the original position plus the number of reinsertions up to its original position minus the number of vehicle failures up to its original position. In constraint set \eqref{c1.10.2}, the final position of a failed vehicle is set to its reinsertion position plus the number of reinsertions up to its reinsertion position minus the number of fails up to its reinsertion position. The last term $y_{vT}\beta_{i}$ is added because the order of the vehicles reinserted to the dummy position $|T|+1$ is irrelevant. 
Constraint set \eqref{c1.20} determines the processing times of vehicles based on the final sequence. We note that the processing time of the dummy reinserted vehicles is set to cycle time for each station (as neutral vehicles). 
Constraint sets \eqref{c1.21} - \eqref{c1.22} determine the starting position and work overload at each position based on the final sequence, respectively. 
Constraint sets \eqref{c1.23} and \eqref{c1.24} guarantee a regenerative production plan. Constraint sets \eqref{c1.25} and \eqref{c1.26} represent the domain of the second-stage variables.

One can obtain the deterministic equivalent formulation (DEF) for MMS with stochastic failures and reinsertion by adding the first-stage constraints \eqref{c1.1.1}-\eqref{c1.3.2} to the second-stage formulation, and by adding copies of all second-stage variables and constraints. 

For our problem, integrating the reinsertion process into MMS presents a bi-objective environment since inserting additional vehicles into an already optimized sequence conflicts with the work overload minimization. The two objectives are given as a single objective in the summation form in \eqref{formulation: reinsertion-second-stage} for the sake of clarity. However, we tackle these two conflicting objectives via a bi-objective approach which is adapted into the proposed solution approaches in Section \ref{section: reinsertion-solution-approach}.
In a multi-objective environment, the objectives conflict with each other, making it impossible to find a solution that optimizes all objectives simultaneously. In order to balance the tradeoffs among the objectives, the concept of \textit{Pareto optimality} is defined. The set of Pareto optimal solutions is called \textit{Pareto front}. A solution is Pareto optimal if the solution is not dominated by any other solutions. The domination criteria is defined as follows. Let $F(x) = f_1(x),..,f_{m}(x)$ be the set of objectives, where $m$ is the number of objectives. Solution $x$ dominates solution $x^{'}$ if $x_{i}\leq x_{i}^{'}$ for all $i\in {1,\ldots,m}$ and $x_{j}\leq x_{j}^{'}$ for at least one of the objectives. The goal of multi-objective algorithms is to find the best representation of the true Pareto front. 

Let us define the objectives separately before proposing our solution approaches.
The \textit{work overload} objective \eqref{eq: reinsertion-first-objective} minimizes total work overload duration for the final sequence (after vehicle failure realization and reinsertions) across all scenarios, stations, and positions.
\begin{flalign}\label{eq: reinsertion-first-objective}
\min \quad  & \sum_{\omega \in \Omega} \rho_{\omega} \sum_{t\in \Bar{T}}\sum_{k\in K} w_{kt\omega}
\end{flalign}
The \textit{reinsertion objective} \eqref{eq: reinsertion-second-objective} minimizes the sum of the squared number of waiting days of the failed vehicles that are not reinserted at the end of the current horizon, across all scenarios. The waiting time is calculated as the number of days between the current production horizon and the horizon that the vehicle was originally planned to be produced. The motivation behind this objective is meeting the delivery deadlines and ensuring customer satisfaction since each car is produced based on a specific order.
\begin{flalign}\label{eq: reinsertion-second-objective}
\min \quad  &  \sum_{\omega \in \Omega}\rho_{\omega} \sum_{i\in F}(1-\sum_{t\in T}y_{it\omega})(g_{i}+1)^2
\end{flalign}

\section{Solution Approaches}\label{section: reinsertion-solution-approach}
MMS, a problem known for its NP-hard nature \citep{tsai1995mixed}, becomes significantly more challenging to solve when faced with stochastic product (car) failures and integrated reinsertion process in a multi-objective setting. Therefore, it is essential to develop effective heuristic methods to tackle large-scale problems.
Accordingly, in this section, we first explain a feasibility enhancement procedure that is utilized in the evolutionary algorithm, and some efficiency improvements that are utilized for all three proposed solution approaches, then we explain our sampling approach that tackles the exponentially increasing number of scenarios. Next, we propose three efficient metaheuristic approaches to solve our problem. First, a population-based evolutionary optimization algorithm utilizing the NSGA-II structure is proposed in Section \ref{section: reinsertion-evolutionary-algorithm}. Then, a two-stage multi-objective local search is presented in Section \ref{section: reinsertion-local-search}. Finally, a hybrid approach, a local search integrated evolutionary algorithm utilizing NSGA-II structure, is proposed in Section \ref{section: reinsertion-local-search-integrated-evolutionary}. Evolutionary and local search algorithms are proposed since they have proved their success on multi-objective problems and MMS problems, respectively. 

\paragraph{Feasibility Enhancement}
Each solution (sequence) is feasible in a standard MMS problem. However, in our problem, with the addition of vehicle failures and the reinsertion process, it is not the case. Although each first-stage decision (sequence before failures and reinsertion) is feasible, the second-stage decisions (final sequence after reinsertions for each scenario) can be infeasible due to one of the following five restrictions regarding the constraints \eqref{c1.2.1} - \eqref{c1.2.5}. The first three restrictions \eqref{c1.2.1} - \eqref{c1.2.3} are tackled while making the decision of the reinsertion position to ensure that a failed vehicle is not reinserted before it is ready, referred to as \textit{being ready feasibility}. On the other hand, the latter two restrictions \eqref{c1.2.4} and \eqref{c1.2.5} are tackled once the final sequence is determined. In case of a restriction \eqref{c1.2.4} violation, which means the number of not reinserted failed vehicles exceeds the limit, we check all positions until we find appropriate positions for a number of vehicles equal to the excess amount. If we cannot accomplish enough reinsertions, then we do the reinsertions with the sacrifice of violating the restriction \eqref{c1.2.5}. Hence, the only restriction that we allow to be violated is \eqref{c1.2.5}, referred to as \textit{lambda feasibility}, however, in case of a violation of this restriction, we try to eliminate (or minimize) it by exploring possible movements (one vehicle movement at a time) without violating other restrictions. In the end, we keep the number of violated subsequences (the subsequences length of $\lambda$ with more than one reinsertion) in order to use the constrained-domination process within the evolutionary algorithms.

\paragraph{Efficiency Improvements}
The work overload and the reinsertion objective function values are calculated as given in \eqref{eq: reinsertion-first-objective} and \eqref{eq: reinsertion-second-objective}, respectively. The computational cost associated with evaluating the objective function and checking and maintaining the feasibility is high. These procedures consume a significant portion of the algorithm's execution time, potentially limiting their overall efficiency and speed in finding a good representation of the Pareto front.
We utilize several improvement techniques to improve the computational efficiency of the proposed algorithms. We employ an accelerated objective evaluation technique to speed up the objective evaluation of each new solution and tabu rules to increase the chance of finding a good reinsertion position, as given in \citet{yilmazlar2023mixed}. That is, every time we generate a random reinsertion position or every time we apply a transformation operator to the second-stage solution during the local search, we apply the tabu rules. Finally, in order to decrease the computational cost of maintaining feasibility, we keep track of the first-stage and reinsertion positions of failed vehicles, and also a reinsertion map of the sequence for each scenario. The reinsertion map ($rm$) is a positional binary mapping of the sequence, $rm_{t}=1$ if there is a reinserted vehicle at position $t$, 0 otherwise. Accordingly, we obtain three improvements as follows: (1) having reinsertion positions of failed vehicles decreases the single iteration cost of the second-stage improvement from $O(N)$ to $O(1)$, (2) having positions of failed vehicles at the first-stage solution decreases the maintaining being ready feasibility cost from $O(N)$ to $O(1)$, (3) having reinsertion position of failed vehicles decreases the maintaining lambda feasibility cost from $O(|F|\lambda)$ to $O(\lambda)$. We note that while improvements (2 and 3) are for all three approaches, improvement (1) is for the local search and local search integrated evolutionary approaches.

\paragraph{Sampling Approach}
Due to the exponential increase in the number of scenarios, the stochastic program becomes overwhelmingly large. Consequently, we utilize a sampling approach to address this challenge. We notice that the majority of failure scenarios have marginal probabilities. Rather than explicitly exploring all potential failure scenarios $\omega \in \Omega$, we approximate the expected value function of the two-stage stochastic program $z^*= \min_{x \in X}  \mathbb{E}[Q(x,\xi_\omega)]$ with an identically and independently distributed (i.i.d) random sample of $N$ realizations of the random vector $\Omega_{N}:=\{\omega_1,\ldots,\omega_N\}\subset \Omega$.
We no longer take into account the probabilities associated with each scenario. However, it is possible for more probable scenarios to occur multiple times within our sample $\Omega_{N}$. We assess the average work overloads and wait times in the concluding sequences of all failure scenarios included in the sample. Let $\hat{N}$ and $n_{\omega}$ represent the set of unique scenarios in $\Omega_N$ and the number of their occurrences. Hence, the objective function \eqref{c1.0.1} $\sum_{\omega \in \Omega} \rho_\omega (\cdot)$ changes to $\frac{1}{N}\sum_{\omega \in \Omega_{N}} (\cdot)$ or equivalently, $\frac{1}{N}\sum_{\omega \in \hat{N}} n_\omega(\cdot)$.

\subsection{Evolutionary Optimization Algorithm}\label{section: reinsertion-evolutionary-algorithm}
We employ the NSGA-II framework which is developed by \citet{deb2002fast}. NSGA-II has demonstrated its efficiency in solving multi-objective problems in various fields \citep{wang2020gaussian, liu2020bi, xu2020data}. 

In order to tackle the second-stage feasibility issue, we employ two procedures. First, we apply \textit{feasibility enhancement procedure} that transforms (if possible) an infeasible sequence into a feasible sequence as explained earlier in this section. Second, we employ the constraint-domination principle defined in \citep{deb2002fast}, which considers the feasibility of each sequence while deciding domination between two solutions. 
Utilizing the constrained-domination principle leads to the outcome where feasible solutions surpass infeasible solutions in terms of nondomination rank. Feasible solutions are evaluated and ranked based on their nondomination level using objective function values. Nevertheless, when comparing two infeasible solutions, the solution with a lower degree of constraint violation achieves a higher rank. In case two infeasible solutions have the same degree of constraint violation, then they are compared based on their objective function values.

In addition to the improvement procedures that we employed for the NSGA-II framework as mentioned above, the novelty of the proposed approach lies in the problem-specific chromosome design and next-generation creation methods: crossover and mutation. 

\paragraph{Initial Population, First-Stage}
We generate the first-stage initial population by employing iterative greedy heuristics. One solution is generated using the greedy heuristic given by \citet{yilmazlar2023mixed}, let us refer to it as the utilization rate greedy heuristic throughout this paper. The rest is generated by using a naive greedy heuristic which makes one assignment at a time, starting from the first position. The vehicle for the first position is selected randomly, then, for the next position, the vehicles that cause the minimum work overload are determined. Out of this group, the vehicles that cause the minimum idle time are determined. A random vehicle is selected from the final group and assigned to the next position. This procedure is repeated until the sequence is completed.

\paragraph{Initial Population, Second-Stage}
The second stage solutions are generated in two steps. We first determine the binary reinsertion decisions of each failed vehicle for each scenario $\omega \in \Omega_{N}$. Then, the reinsertion positions of the failed vehicles that are decided to be reinserted are determined. Note that if a failed vehicle is decided to be reinserted but it does not become ready before the end of the horizon, then the binary reinsertion decision is changed. The feasibility enhancement procedure is applied to each scenario of each solution, once the second-stage decisions are determined. All second-stage decisions are randomly determined (both binary reinsertion and reinsertion position decisions), except for two solutions: 1. none of the failed vehicles are reinserted (except the ones that have to be reinserted) for the solution that is generated using the  utilization rate greedy heuristic, 2. all of the failed vehicles are reinserted (except the ones that cannot be reinserted due to infeasibility) for one of the other solutions. 

\begin{figure}[h]
\captionsetup{justification=centering}
    \centerline{\includegraphics[scale=0.8]{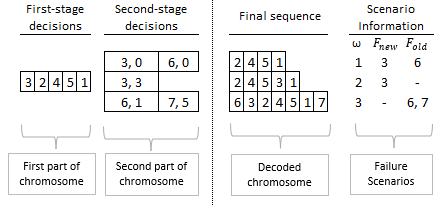}}
    \caption{Chromosome illustration over an example with three scenarios}\label{fgr: reinsertion-chromosome-design}
\end{figure}

\paragraph{Chromosome Design}
The chromosome representation of each solution is designed in two parts which represent the first and second-stage decisions, respectively. As illustrated in Figure \ref{fgr: reinsertion-chromosome-design}, while the first part consists of a single sequence, the second part consists of $|\Omega_{N}|$ sequences, each representing a final sequence for scenario $\omega \in \Omega_{N}$, note that we illustrate only three scenarios in the figure.
The first part of the chromosome encodes the corresponding first-stage sequence. Then, vehicle failures are realized which is summarized under the scenario information column. Next, second-stage decisions (reinsertions) are encoded in the second part of the chromosome. From the scenario information, we know that vehicles $\{3,6\}$, $\{3\}$, and $\{6,7\}$ are in $F$ for the scenarios 1,2, and 3, respectively. The first and second numbers of each gene (tuple) correspond to the id of a failed vehicle and the reinsertion decision for that vehicle, respectively. If a failed vehicle is not reinserted, then the second number equals zero, otherwise equals the reinsertion position. For example, in scenario 2, the failed vehicle 3 is reinserted in position 3. 

\begin{figure}[h]
\captionsetup{justification=centering}
\centering
\begin{subfigure}{0.4\textwidth}
  \centering
    \includegraphics[width=0.6\textwidth]{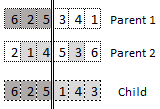}\hfill
      \caption{PMX for the first part}
  \label{fig_sub: reinsertion-pmx}
\end{subfigure}%
\begin{subfigure}{0.6\textwidth}
  \centering
    \includegraphics[width=0.35\textwidth]{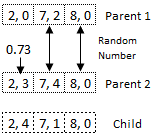}\hfill
  \caption{Uniform crossover for the second part of a given scenario}
  \label{fig_sub: reinsertion-uniform-crossover}
\end{subfigure}%
\caption{Illustration of proposed crossover methods}
\label{fig: reinsertion-crossovers}
\end{figure}

\paragraph{Crossover}
The crossover is executed for the first and second parts of the chromosomes separately. First, two solutions from the parent population are randomly selected. Then, in order to generate a single child by using the two selected parents, a modified partially mapped crossover (PMX) is applied to the first part of the chromosome as illustrated in Figure \ref{fig_sub: reinsertion-pmx}. In order to generate a single child chromosome, a random single point is selected. While the genes up to that point are selected from the first parent, the rest of the genes are selected from the second parent, based on the order that occurred on the second parent, and the genes were not selected from the first parent. 

The second part of a child chromosome is generated in two stages. First, a uniform crossover is employed to decide the decimal numbers in the second part of the child's chromosome, however, as a binary selection on binary reinsertion decision. Thus, for each gene in each scenario, if the binary reinsertion decisions are the same in two-parent genes, the child gets the same decision. Otherwise, we generate a random number to decide on the dominant parent for that specific gene. Next, the reinsertion position is determined independently from parents since otherwise, most of the children chromosomes become infeasible. 
The second part crossover procedure is illustrated in Figure \ref{fig_sub: reinsertion-uniform-crossover}. The parents have the same binary decision for the second and third genes, thus the child gets the same binary reinsertion decision, vehicle 7 is reinserted and vehicle 8 is not reinserted. The parents have different decisions for the first gene, thus a random number is generated to decide which decision the child gets. Since the random number is greater than $0.5$, the child gets the decision from the second parent which is reinserting vehicle 2. The reinsertion positions of vehicles 2 and 7 are randomly decided independent of corresponding parent genes, but obeying being ready feasibility restrictions. 

\paragraph{Mutation}
We apply mutations to the first and second-stage solutions separately. For the first-stage solutions, an inversion operator is employed as a mutation in order to avoid local optima and to increase the diversity \citep{akgunduz2010adaptive}. The inversion operator is applied to a randomly selected subsequence (the subsequence between randomly selected two points). 
For the second-stage solutions. a single gene mutation is applied to a gene that carries the binary reinsertion decision. Accordingly, we generate a random number for each child chromosome, if the random number is smaller than the mutation threshold, then we apply mutation to a randomly selected gene on a randomly selected scenario. On one hand, if the selected gene is a negative binary reinsertion decision (not reinserted), then it changes to a positive decision (reinsert), so that the corresponding failed vehicle is reinserted to a random position. On the other hand, if the selected gene is a positive binary reinsertion decision, then it changes to a negative decision, so that the corresponding vehicle is removed from the final sequence and the succeeding vehicles are moved forward to close the gap.

\subsection{Two-stage Bi-objective Local Search Algorithm}\label{section: reinsertion-local-search}
This section proposes a simulation-based two-stage bi-objective local search (STMLS) algorithm. A pseudo-code for STMLS is given in Algorithm \ref{alg: reinsertion-pseudo-code-stmls}. 
To summarize the algorithm, STMLS gets an initial solution that is generated for the deterministic counterpart of the proposed MMS problem, which excludes vehicle failures and reinsertions, this problem is referred to as \textit{one-scenario problem}. Then, the initial solution is improved over the one-scenario problem for just a couple of seconds. Next, the improvement loop for the full-information problem starts once the second-stage solutions are generated. During this loop until meeting the termination criterion, STMLS alternates between improving the first-stage and second-stage solutions on work overload objective and introduces random changes to the binary reinsertion decisions to explore the reinsertion objective space. The algorithm aims to find a good representation of the Pareto front by solving the proposed problem in a bi-objective environment.

\begin{algorithm}[H]
\scriptsize
\caption{Pseudo-code of STMLS}\label{alg: reinsertion-pseudo-code-stmls}
\begin{algorithmic}
\STATE{\textbf{Input:} The initial first-stage solution generated by utilization rate greedy heuristic, parameters: operator weights, $\theta$,  $\tau_{f}$, and $\tau_{s}$.}
\STATE{\textbf{Output:} External population: a set of non-dominated solutions that represent the Pareto front}
\STATE{Improve the initial first-stage solution over one-scenario problem}
\STATE{External population $\leftarrow \emptyset$}
\STATE{Generate random second-stage solutions for each $\omega \in \Omega_{N}$ (all possible binary reinsertion decisions are 1), and apply feasibility enhancement}
\STATE{Calculate the objective function value of each solution in the initial population}
\STATE{Improve second-stage solutions on work overload objective by applying second-stage solution improvement procedure}
\STATE{iteration $\leftarrow 1$}
\WHILE{termination criterion is not reached}
    \IF{iteration \% $\theta = 0$}
        \STATE{Update external population}
        \FOR{each $\omega \in \Omega$}
            \STATE{Select a failed vehicle randomly and change the binary reinsertion decision}
            \STATE{If the change is from 0 zero to 1, then find a random feasible position for the reinsertion (if possible)}
            \STATE{Calculate the new objective function values and accept the new solution as the incumbent solution, regardless of objective values}
            \STATE{Apply second-stage improvement procedure for $\tau_{s}$ duration}
        \ENDFOR
    \ELSE
        \STATE{Apply first-stage improvement procedure for $\tau_{f}$ duration}
    \ENDIF
    \STATE{iteration $\leftarrow iteration + 1$}
\ENDWHILE
\end{algorithmic}
\end{algorithm}

\paragraph{Local Search Specific Efficiency Improvements}
Since the movements made via transformation operators make local changes to the sequence, there is no need to reevaluate the whole sequence. Hence, we utilized partial objective reevaluation during a local search (both first-stage and second-stage), in addition to the accelerated reevaluation, as given in \citep{yilmazlar2023mixed}. 

\paragraph{Transformation Operators}
We apply one randomly selected operator at each iteration of the improvement procedures explained below. The random selection of operators is based on the operator weights which are determined based on our preliminary experiments. The swap operator is used to interchange the positions of two randomly chosen cars in the sequence. The insertion operator removes a car from its current position $i$ and inserts the car to a position $j$. There are two types of insertions: backward insertion and forward insertion. Backward insertion occurs when $i>j$, where the car is inserted at position $j$ and all vehicles between positions $j$ and $i$ are shifted one position to the right (scheduled later). On the other hand, forward insertion happens when $i<j$, where the car is inserted at position $j$ and all vehicles between positions $i$ and $j$ are shifted one position to the left (scheduled earlier). Inversion reverses a subsequence between two randomly selected positions in the whole sequence. Even though there are more operators proposed in the literature that could be used to transform a permutation-based representation, a local search across the neighborhood defined over these four transformation operators is shown to be very efficient to solve sequencing problems \citep{estellon2006large, yilmazlar2023adaptive}. 

\paragraph{First-stage Improvement Procedure}
Given a full solution (first and second-stage solutions), this procedure executes a local search within the neighborhood defined above, by applying one operator at a time on the first-stage solution. Each time a neighbor first-stage solution is visited, the feasibility of each second-stage solution is checked and a new reinsertion position is randomly generated (if possible) for the ones that cannot be reinserted to the previously determined position anymore due to the change on the first-stage solution. We note that tabu rules are utilized while selecting the reinsertion positions during the random second-stage solution generation process, additionally, the binary reinsertion decisions are not changed during this procedure unless due to a feasibility issue. If the new solution is not deteriorated based on the work overload objective, then the new solution is accepted as the incumbent solution, rejected otherwise. To summarize, this procedure improves the work overload objective while keeping the reinsertion objective unchanged (when possible). 

\paragraph{Second-stage Improvement Procedure}
Given a second-stage solution for a scenario $\omega \in \Omega_{N}$, this procedure executes a local search on the second-stage reinsertion position decisions by applying swap and insertion operators on the reinserted vehicles only. That is because, a movement of vehicle $v\in V$ means changing the first-stage solution, and there is no valid inversion movement when $\lambda>1$. A new solution is accepted if it is not deteriorated in terms of work overload objective. 

\subsection{Local Search Integrated Evolutionary Algorithm}\label{section: reinsertion-local-search-integrated-evolutionary}
In this section, we propose a hybrid approach which is a local search integrated evolutionary algorithm, \textit{LS-NSGA-II}. We integrate local search improvement procedures to the NSGA-II structure proposed in Section \ref{section: reinsertion-evolutionary-algorithm} in two places. First, each first-stage solution in the initial population is improved over the one-scenario problem. Second, the second-stage improvement procedure is applied to each child chromosome when it is generated. Executing a local search on each child chromosome increases the computational complexity of the algorithm, hence, it is necessary to balance the time spend on local search and evolution. In order to accomplish this, we decrease the population size $P$ and set a shorter time limit $\tau_{s}$ for the second-stage improvement on each chromosome. A pseudo-code for the proposed approach is given in Algorithm \ref{alg: reinsertion-pseudo-code-nsga}.

\begin{algorithm}[H]
\scriptsize
\caption{Pseudo-code of Local Search Integrated NSGA-II}\label{alg: reinsertion-pseudo-code-nsga}
\begin{algorithmic}
\STATE{\textbf{Input:} The first stage initial population generated by greedy heuristic, parameters: population size ($P$), mutation probability, and $\tau_{s}$}
\STATE{\textbf{Output:} A set of non-dominated solutions that represent the Pareto front}
\STATE{*Improve each first-stage solution over one-scenario problem}
\STATE{*EP $\leftarrow \emptyset$}
\STATE{For each first-stage solution, randomly generate second-stage solutions for each $\omega \in \Omega$, and apply feasibility enhancement}
\STATE{Calculate the objective function value of each solution in the initial population}
\WHILE{Termination criterion is not reached}
    \STATE{Create a child population of size $P$ using the parent population via applying crossover, mutation, and feasibility enhancement}
    \STATE{*Apply second stage improvement procedure, given in Section \ref{section: reinsertion-local-search}, to each child chromosome for $\tau_{s}$ duration}
    \STATE{Calculate the objective function value of each child solution}
    \STATE{Combine parent and child populations which results in a set of solutions in size $2P$}
    \STATE{Determine nondomination rankings using the constrained domination principle}
    \STATE{Calculate crowding-distance of each solution}
    \STATE{Select the best $P$ solutions (next generation's parent population) based on the nondomination ranking and crowding distance in the given order}
    \STATE{*Update EP: Add all non-dominated solutions to the EP and remove all dominated solutions from the EP}
\ENDWHILE
\end{algorithmic}
\footnotesize{\emph{*Only with the LS integrated NSGA-II structure}}
\end{algorithm}

One of the drawbacks of NSGA-II is that the number of non-dominated solutions (solutions in the first Pareto front) is limited with the population size $P$. In order to tackle this drawback we collect all non-dominated solutions that are explored along the search in an external set (external population, EP), which is originally proposed by \cite{michalak2015improving}. We note that, we utilize this improvement only with the LS-NSGA-II structure since our preliminary experiments show that the number of non-dominated solutions does not exceed the population size with the standard NSGA-II structure.

\section{Numerical Experiments}\label{section: reinsertion-numerical-experiments}
In this section, we first describe the experimental setup in Section \ref{section: reinsertion-experimental-setup}. Then, the computational performance of the proposed solution approaches is assessed in Section \ref{section: reinsertion-computational-performance}. Next, the solution quality of the solutions obtained by solving the proposed problem is assessed in Section \ref{section: reinsertion-solution-quality}, over dynamic reinsertion simulations, by comparing them with the solutions obtained by solving the one-scenario problem and the problem given in \citet{yilmazlar2023mixed}.

We note that the traditional MMS problem is NP-hard \cite{tsai1995mixed}, hence, the problem proposed in this study is NP-hard since it is a more complicated version of the traditional MMS. Considering the stochastic failure of vehicles poses a significant complexity increase in the model. Additionally, integrating vehicle reinsertions to MMS adds significant complexity. Accordingly, we did not propose any exact solution approaches since they are inadequate to solve even small instances. e.g., Gurobi can optimally solve instances with up to 5 vehicles, 5 stations, and 100 scenarios even with a single objective, in a one-hour time limit.

\subsection{Experimental Setup}\label{section: reinsertion-experimental-setup}
We generated real-world inspired instances from our automobile manufacturer partner's assembly line and planning information, the numbers of vehicles are 200, 300, and 400. We generated 30 instances for each parameter configuration, 90 instances in total. All instances have five stations, of which the first station is selected as the most restrictive station for EVs, the battery loading station. The remaining four stations are selected among other critical stations that conflict with the battery loading station. 

The cycle time $c$ is 97 TU, and the station length $l$ is 120 TU for all but the battery loading station, which is two station lengths, 240 TU. Table \ref{tbl:mms-processing-times} provides information about the distribution of the processing times. The average and maximum processing times for each station are lower than the cycle time and the station length, respectively. Additionally, the ratio of the EVs is in the range of [0.25, 0.33] among all vehicles, for all instances. 

\begin{table}[h]
\captionsetup{format=plain, labelfont=bf, justification=centering, singlelinecheck=false, font=small, skip=4pt}
\caption{Processing times distribution}
\centering
\scalebox{0.8}{
\begin{tabular}{
>{\columncolor[HTML]{FFFFFF}}c 
>{\columncolor[HTML]{FFFFFF}}c 
>{\columncolor[HTML]{FFFFFF}}c 
>{\columncolor[HTML]{FFFFFF}}c}
\multicolumn{1}{l}{\cellcolor[HTML]{FFFFFF}}           & \multicolumn{3}{c}{\cellcolor[HTML]{FFFFFF}Time (s)} \\ \cline{2-4} 
\multicolumn{1}{l}{\cellcolor[HTML]{FFFFFF}Station ID} & Min             & Mean           & Max  \\  \hline  
1                                                      & 42.6            & 94.1                          & 117.2           \\
2                                                      & 7.9             & 84.3                          & 197.9           \\
3                                                      & 57.8            & 96.2                          & 113.3           \\
4                                                      & 26.9            & 96.9                          & 109.7           \\
5                                                      & 57.8            & 96.2                          & 114.3          
\end{tabular}}
\label{tbl:mms-processing-times}
\end{table}

The maximum number of failed cars allowed at the end of the planning horizon is set to $f_{max}=|V|*0.05$. The reinsertion window length $\lambda=10$. The order lead time is 9. The number of days between the planned and current period for the failed car $i\in F_{old}$, $g_i$, is discrete uniformly distributed in [1, 9]. The number of extra days until the delivery date that the failed car $i\in F_{old}$, $d_i$, is discrete uniformly distributed in [$g_{i}$, 9]. The number of positions required for the failed vehicle $i$ to be ready for reinsertion, $r_i$, is discrete uniformly distributed in [10, $|V|$-10] for $i\in F_{new}$ and in [0, $|V|$-10] for $i\in F_{old}$. As aforementioned, these parameters are generated for each instance as a deterministic parameter, in other words, they are the same for each scenario $\omega \in \Omega$.

The failure rates were derived from six months of historical data by performing predictive feature analysis on vehicles. Based on the analysis, two groups of vehicles are formed according to their failure probabilities, low-risk and high-risk vehicles, whose failure probabilities are in the range of [0.0, 0.01] and [0.2, 0.35], respectively. The ratio of high-risk vehicles was set to [0.03, 0.05] among all vehicles. We utilized a sampling approach by generating an i.i.d random sample of $N$ realizations of failure scenarios. For each failure scenario $\omega \in \Omega_N$ and vehicle $v\in V$, we first chose whether the vehicle was high-risk or low-risk (based on their prevalence). Then, depending on being a high-risk or low-risk vehicle, a failure probability was randomly selected from the respective range. Next, it was determined whether the vehicle failed or not. We did not allow any low-risk vehicle to fail at any scenario in order to have a more representative sample of scenarios. Finally, the number of vehicles in $F_{old}$ is determined for each scenario from a discrete uniform distribution in [0, $f_{max}$]. Once $|F_{old}|$ is determined for each scenario, the vehicles in $F_{old}$ are randomly selected from a predetermined set of vehicles in size $f_{max}$. We set the sample size $N=100$ for all experiments unless otherwise stated, considering the problem complexity and the results given in \citep{yilmazlar2023mixed}.  

The proposed algorithms are implemented in Python 3. In order to meet the industry requirements that expect to obtain a good quality solution in a limited time, the termination criterion is set to 600 seconds for all of the algorithms. The number of replications (runs) is set to 30, i.e., each instance is solved by each solution approach for 30 times with a different seed. The computing nodes, with 16 cores and 125 GB of memory, of the Clemson University supercomputer are utilized to execute the numerical experiments. 

\subsection{Computational Performance}\label{section: reinsertion-computational-performance}
In this section, we assessed the performance of the proposed algorithms: STMLS, NSGA-II, and LS-NSGA-II. The performance assessment is executed in two parts. We first evaluate the solution approach on each objective, separately. Then, a multi-objective comparison over the non-dominated sets provided by each algorithm is executed.

\begin{figure}[h]
\captionsetup{justification=centering}
\centering
\begin{subfigure}{0.5\textwidth}
  \centering
  \includegraphics[width=0.85\linewidth]{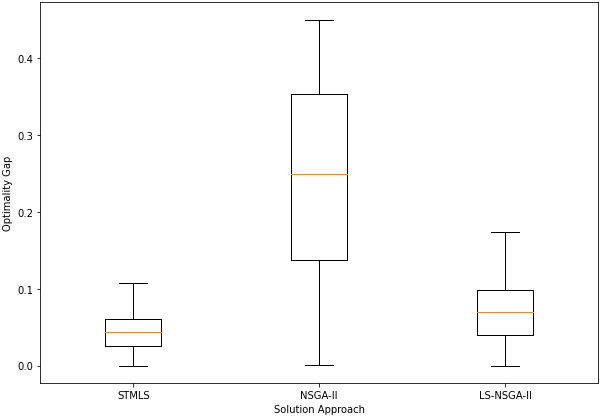}
  \caption{work overload objective}
  \label{fig_sub: reinsertion-wo-objective-convergence}
\end{subfigure}%
\begin{subfigure}{0.5\textwidth}
  \centering
  \includegraphics[width=0.85\linewidth]{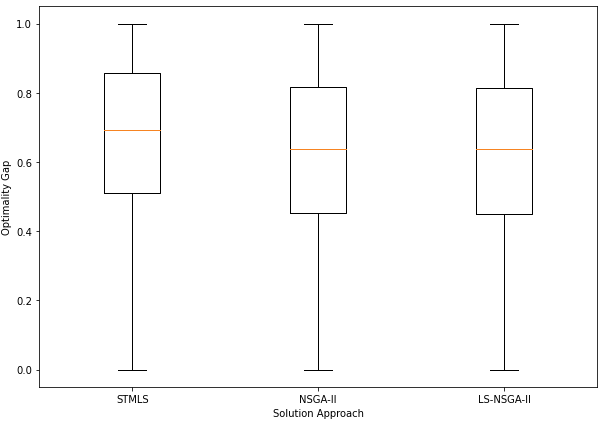}
  \caption{reinsertion objective}
  \label{fig_sub: reinsertion-wait-objective-convergence}
\end{subfigure}
\caption{Convergence comparison of the algorithms in terms of each objective separately}
\label{fig: reinsertion-convergence-of-algorithms-single-objective}
\end{figure}

\paragraph{Single-objective Analysis}
After completing each run of every solution method, we obtained a set of non-dominated solutions: the first Pareto front of the last generation for the NSGA-II and the external populations for the STMLS and the LS-NSGA-II. In order to evaluate the reliability of the solution methods for each objective, we then compared each non-dominated set with the best-found value of each objective across all runs. Let us call the best found values across all runs as \textit{heuristic ideal points}, because finding ideal points (the optimal values of each objective when the problem solved explicitly for each objective) is not possible. This comparison is executed as follows: for each instance, we first determine the heuristic ideal points. Then, we determine the best-found objective values of each non-dominated set, one value per objective from each run. Next, the gap between the heuristic ideal points and the best-found in a single run is calculated for each objective. The box plots that show the distribution of the gaps are given in Figure \ref{fig: reinsertion-convergence-of-algorithms-single-objective}. 

For the work overload objective, STMLS and LS-NSGA-II provide more reliable solutions, each within an 11\% and 18\% gap, respectively. While STMLS provides slightly better solutions compared to LS-NSGA-II, there is no significant difference based on the Kruskal-Wallis test with 95\% confidence interval. On the other hand, they both outperform the NSGA-II which provides solutions in a wider objective value range, with up to 44\% gap.
Moreover, in terms of reinsertion objective, all three solution approaches provide similar solution quality based on this analysis. Since, for most of the instances, the best-found reinsertion objective is either zero or very close to zero, even a small fraction of the increase in the reinsertion objective value induces a large gap, as in the figure, each approach have on average more than 60\% gap. However, the average ideal point of reinsertion objective is 0.56 and the average best-found of each non-dominated set is 1.55, 1.34, and 1.34, for the STMLS, NSGA-II, and LS-NSGA-II, respectively. 
We can conclude the single-objective analysis by stating that STMLS and LS-NSGA-II can provide reliable solutions in terms of both objectives and NSGA-II can provide reliable solutions only for the reinsertion objective.

\paragraph{Multi-objective Analysis}
The performance of the proposed solution approaches was compared in a multi-objective optimization environment. That is, a comparison of non-dominated sets of solutions provided by each approach (Pareto front). The quality of the Pareto front generated by an algorithm is assessed based on criteria such as the number of non-dominated solutions explored, the distribution of the solutions in the non-dominated set, and the distance among the non-dominated solutions and Pareto-optimal solutions or ideal points \citep{zitzler2000comparison, alaghebandha2017optimizing, schaffer1985multiple, hyun1998genetic, coello2007evolutionary, srinath2022study}. In this section, we first present empirical attainment surface graphs in order to have an insight into the Pareto exploration and quality comparison of the algorithms. The attainment function provides a description of a random non-dominated point set's location distribution. The attainment function can be estimated by empirical attainment functions (EAF) \citep{fonseca1996performance}, which is utilized to interpret the statistical performance of stochastic multi-objective optimizers. Next, we employed the following performance metrics to assess the multi-objective performance comparison of the solution approaches. We note that the objective function values are normalized for each instance separately and the normalized values are used throughout this section in order to tackle different scales associated with each objective.
\begin{itemize}
    \item \textit{The number of non-dominated solutions (NNS)} evaluates the number of solutions within the Pareto front that is generated by an algorithm, i.e., the number of non-dominated solutions explored by an algorithm.
    \item \textit{Covered size of space (CSS)} compares two sets of non-dominated solutions provided by two algorithms. The CSS metric, also known as \textit{C metric}, is introduced by \citet{zitzler2000comparison}. Given two sets of non-dominated solutions, $X$ and $Y$, $C(X, Y)$ calculates the ratio between the number of solutions in set $Y$ that are dominated by any solution from set $X$ and the total number of solutions in set $Y$, given in \eqref{eq: reinsertion-css}. That is, $C(X, Y)=1$ means that all the solutions in set $Y$ are dominated by at least one of the solutions from set $X$.
    \begin{equation}\label{eq: reinsertion-css}
    C(X, Y)= \frac{\mid a_{Y}\in Y \, ; \, \exists a_{X}\in X, \, a_{X}\leq a_Y \mid}{Y}
    \end{equation}
    \item \textit{Mean ideal distance (MID)} measures the distance between the ideal point, denoted as $z^*$, and a given set of non-dominated solutions $X$, as given in \eqref{eq: reinsertion-mid}, assume that the $m$ is the number of objectives and $n$ is the number of solutions in set $X$. Remember that we employed the best-found objective function values across all runs as heuristic ideal points.
    \begin{equation}\label{eq: reinsertion-mid}
    MID = \frac{\sum_{i=1}^n\sqrt{\sum_{j\in m}(X_{ij}-z_m^*)^2}}{n}
    \end{equation}
    \item \textit{Spread of non-dominance solutions (SNS)} evaluates the distribution of non-dominated solutions based on the MID value. The SNS is used to assess the variability of the distance between each solution and the ideal point. Note that a higher SNS value is desirable as it indicates a more uniform distribution of solutions across the Pareto front.
    \begin{equation}\label{eq: reinsertion-sns}
    SNS = \sqrt{\frac{\sum_{i=1}^n \Big(MID - \sqrt{\sum_{j\in m}(X_{ij}-z_m^*)^2}\Big)^2}{(n-1)}}
    \end{equation}
\end{itemize}

\begin{figure}[h]
\captionsetup{justification=centering}
    \centerline{\includegraphics[width=1.0\textwidth]{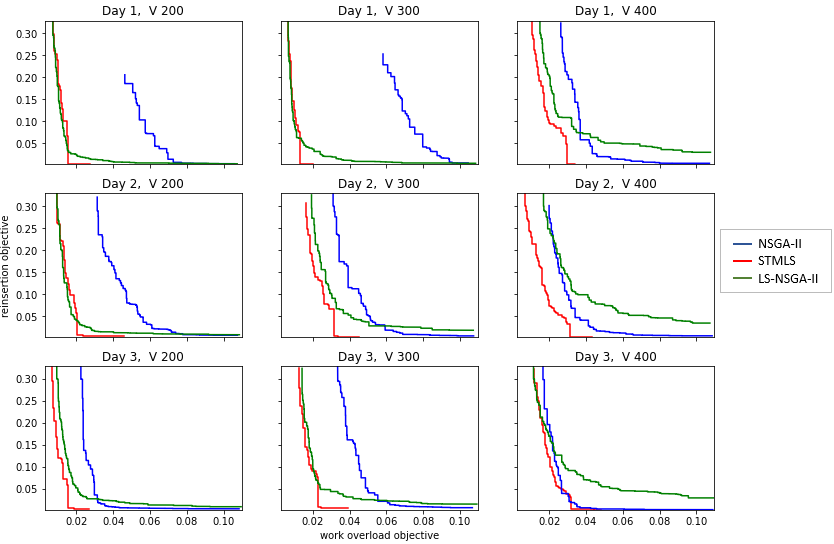}}
    \caption{Comparison of the 50\%-attainment surfaces of the proposed algorithms, with a 95\% confidence level. The illustrated region is limited with up to 0.1 and 0.3 for the work overload and reinsertion objectives, respectively}\label{fgr: reinsertion-attainment-local}
\end{figure}

In Figure \ref{fgr: reinsertion-attainment-local}, we compare the statistical performance of the algorithms over nine instances by utilizing the attainment surfaces. We randomly selected three production days (three instances for each day, each with a different number of vehicles) for this analysis since it is not convenient to share plots for each of the 90 instances. 

The STMLS and LS-NSGA-II algorithms provide better-quality Pareto fronts compared to NSGA-II. It is interesting to note that integrating local search procedures within the NSGA-II structure have improved the NSGA-II algorithm, this was at the expense of degradation in the right-below area, where the work overload objective is high and the reinsertion objective is low, especially as the number of vehicles increases. This is interpreted as the increasing algorithm complexity may hurt the exploration in some areas as the size of the problems increase. 
It is also important to point out that LS-NSGA-II explores a greater region in its Pareto front compared to STMLS and NSGA-II, i.e., better exploration of the edge cases where one of the objectives is very low while the other objective is very high. It is difficult to demonstrate the statistical comparison of the algorithms for whole space in a single graph, hence, Figure \ref{fgr: reinsertion-attainment-local} is limited to focus on the left-below region where more balanced solutions, in terms of objective function values, are provided. Refer to Appendix \ref{appendix: reinsertion-attainment-surface} for the global counterpart of Figure \ref{fgr: reinsertion-attainment-local}.

\begin{table}[h]
\captionsetup{format=plain, labelfont=bf, justification=centering, singlelinecheck=false, font=small, skip=4pt}
\caption{NNS, MID, and SNS performance metrics comparison of the NSGA-II, STMLS, and LS-NSGA-II, averaged across all instances and all runs}
\centering
\scalebox{1.0}{
\begin{tabular}{
>{\columncolor[HTML]{FFFFFF}}l |
>{\columncolor[HTML]{FFFFFF}}c 
>{\columncolor[HTML]{FFFFFF}}c 
>{\columncolor[HTML]{FFFFFF}}c l}
Solution   Method & NNS   & MID  & SNS  &  \\ \cline{1-4}
NSGA-II           & 46.18 & 0.23 & 0.10 &  \\
STMLS             & 12.99 & 0.26 & 0.08 &  \\
LS-NSGA-II        & 40.76 & 0.24 & 0.12 & 
\end{tabular}}
\label{tbl: reinsertion-performance-metrics}
\end{table}

In Tables \ref{tbl: reinsertion-performance-metrics} and \ref{tbl: reinsertion-CSS-metric}, we compare the algorithms in terms of the performance metrics that are explained earlier in this section. The presented values correspond to the average values across all runs and instances, yet for the CSS metrics, the instances are aggregated based on instance size since the performance of the algorithms, for the CSS metric, is impacted by the size of the instances. Table \ref{tbl: reinsertion-performance-metrics} shows that the STMLS performs worse in terms of all three metrics, NNS, MID, and SNS, compared to the evolutionary algorithms. STMLS generates less number of non-dominated solutions. While the generated solutions have a greater distance to the heuristic ideal points, the smaller SNS value shows that the spread of the solutions is low. This is because the non-dominated solutions provided by STMLS have, on average, high reinsertion objective value since STMLS spends most of the computational effort on the minimization of the work overload objective. On the other hand, LS-NSGA-II and NSGA-II perform similarly in terms of NNS, MID, and SNS metrics, yet higher SNS and lower MID values of LS-NSGA-II support that it explores more solutions in the edge regions.

\begin{table}[h]
\captionsetup{format=plain, labelfont=bf, justification=centering, singlelinecheck=false, font=small, skip=4pt}
\caption{CSS performance metric comparison of the NSGA-II, STMLS, and LS-NSGA-II, averaged across all instances and all runs}
\centering
\scalebox{1.0}{
\begin{tabular}{ll|ccc}
Solution   Method & V   & NSGA-II & STMLS & LS-NSGA-II \\ \hline
NSGA-II           & 200 & -       & 0.07  & 0.22       \\
                  & 300 & -       & 0.06  & 0.23       \\
                  & 400 & -       & 0.05  & 0.27       \\ \hdashline
STMLS             & 200 & 0.84    & -     & 0.61       \\
                  & 300 & 0.86    & -     & 0.65       \\
                  & 400 & 0.88    & -     & 0.80       \\ \hdashline
LS-NSGA-II        & 200 & 0.72    & 0.31  & -          \\
                  & 300 & 0.71    & 0.26  & -          \\
                  & 400 & 0.66    & 0.12  & -         
\end{tabular}}
\label{tbl: reinsertion-CSS-metric}
\end{table}

CSS metric provides a reliable pair-wise comparison among sets of non-dominated solutions. Table \ref{tbl: reinsertion-CSS-metric} shows that STMLS consistently achieves high CSS values, which means that STMLS provides better quality Pareto fronts. Although STMLS provides a significantly less number of non-dominated solutions compared to evolutionary algorithms, the high CSS values show that STMLS-provided solutions dominate, on average, 85\% of NSGA-II-provided solutions and 68\% of LS-NSGA-II-provided solutions. LS-NSGA-II performs relatively well compared to NSGA-II.  LS-NSGA-II-provided solutions dominate, on average, 70\% of NSGA-II-provided solutions, so we can say that integration of local search improvement procedure into the NSGA-II structure drastically enhanced the algorithm performance, yet with some decrease with the increasing instance size.

\subsection{Solution Quality}\label{section: reinsertion-solution-quality}
In this section, we evaluated the solution quality of the solutions obtained by solving the problem proposed in this study, \textit{full-information problem with failures and reinsertion} (FFR), over dynamic reinsertion simulations by comparing them with the solutions obtained by solving the one-scenario problem and the problem given in \citep{yilmazlar2023mixed}, \textit{full-information problem with failures} (FF). We utilized the first part of the TS algorithm given in \citep{yilmazlar2023mixed} to solve the one-scenario problem. The complete version of the TS algorithm is employed to solve the FF problem. Finally, the STMLS algorithm is used to solve the FFR problem. We preferred  STMLS over LS-NSGA-II due to the reliable performance of STMLS in terms of work overload objective.

The dynamic reinsertion simulations are executed based on the industry application which is to reinsert a reinstating vehicle as soon as a convenient position is found within the sequence. The simulation procedure is explained as follows: we first solved each problem with the corresponding solution method for all 90 instances and for 30 runs, within a time limit of 600 seconds, and with a sample size $N=100$, note that the sample size is valid only for the full-information problems. Then, another set of 1000 samples is generated as test simulations, and the solutions obtained from each method are simulated over these samples by executing reinsertions dynamically. The reinsertions are made in different work overload thresholds, while also obeying feasibility rules defined in this study. For each solution, each position is checked for each reinstating vehicle, starting from the first position. Once a convenient position for a reinstating vehicle is found that causes less than or equal to the work overload threshold, reinsertion is made. The simulation for a solution is completed when all the positions are checked and the possible reinsertions are made. 

\begin{table}[h]
\captionsetup{format=plain, labelfont=bf, justification=centering, singlelinecheck=false, font=small, skip=4pt}
\caption{Comparison of variants of MMS problem over dynamic reinsertion simulations}
\centering
\scalebox{1.0}{
\begin{tabular}{
>{\columncolor[HTML]{FFFFFF}}c |
>{\columncolor[HTML]{FFFFFF}}c 
>{\columncolor[HTML]{FFFFFF}}c |
>{\columncolor[HTML]{FFFFFF}}c 
>{\columncolor[HTML]{FFFFFF}}c |
>{\columncolor[HTML]{FFFFFF}}c 
>{\columncolor[HTML]{FFFFFF}}c }
\multicolumn{1}{l|}{\cellcolor[HTML]{FFFFFF}}                & \multicolumn{2}{c|}{\cellcolor[HTML]{FFFFFF}One-scenario} & \multicolumn{2}{c|}{\cellcolor[HTML]{FFFFFF}\begin{tabular}[c]{@{}c@{}}Full-information \\      with failures\end{tabular}} & \multicolumn{2}{c}{\cellcolor[HTML]{FFFFFF}\begin{tabular}[c]{@{}c@{}}Full-information\\      with failures and \\      reinsertion\end{tabular}} \\
\begin{tabular}[c]{@{}c@{}}WO \\      Threshold\end{tabular} & Obj WO                      & Obj RE                      & Obj WO                                                              & Obj RE                                                     & Obj WO                                                                  & Obj RE                                                                  \\ \hline
0                                                            & 214.96                      & 113.96                      & \textbf{178.94}                                                     & 133.15                                                     & 187.18                                                                  & \textbf{107.15}                                                         \\
3                                                            & 225.38                      & 90.32                       & \textbf{186.28}                                                     & 106.39                                                     & 191.23                                                                  & \textbf{82.33}                                                          \\
5                                                            & 232.55                      & 81.91                       & \textbf{191.75}                                                     & 97.16                                                      & 195.65                                                                  & \textbf{73.51}                                                          \\
10                                                           & 255.43                      & 50.51                       & 212.27                                                              & 60.87                                                      & \textbf{200.48}                                                         & \textbf{41.36}                                                          \\
15                                                           & 282.07                      & 19.29                       & 239.47                                                              & 23.80                                                      & \textbf{212.45}                                                         & \textbf{16.32}                                                          \\
30                                                           & 317.15                      & 1.84                        & 276.34                                                              & 1.82                                                       & \textbf{238.15}                                                         & \textbf{1.52}                                                          
\end{tabular}}
\label{tbl: reinsertion-simulation-results}
\end{table}

Table \ref{tbl: reinsertion-simulation-results} presents the numerical results of the solution quality simulations. The objective function values were averaged across all instances and all runs and for the FFR problem across all solutions in the corresponding non-dominated set. The \textit{Obj WO} and \textit{Obj RE} columns correspond to the work overload and reinsertion objectives, respectively. The results show that the one-scenario problem solutions are outperformed, in terms of work overload objective, by the FF problem by around 17\% over the simulations with the work overload threshold of 0, 3, and 5, and by the FFR problem by 22\%, 22\%, and 25\% when the work overload threshold is 10, 15, 30, respectively. It is important to point out that the FF problem solutions outperform the FFR problem solutions, in terms of work overload objective, by around 4\% when the work overload threshold is less than 10. It aligns with the problem characteristics. FF problem solutions provide better work overload results when the work overload threshold is low, the reinsertion is not made at the cost of a large work overload sacrifice. It is because the FF problem minimizes work overload over scenarios that consider the failures but not the reinsertions. On the other hand, the FFR problem solutions outperform the FF problem solutions, in terms of work overload, by 6\%, 11\%, and 14\% when the work overload threshold is 10, 15, and 30, respectively. Additionally, FFR problem solutions provide the best reinsertion objective results across all work overload thresholds. Accordingly, we can conclude that the FFR problem can generate good-quality solutions in terms of both objectives, since the motivation of the FFR problem is to generate robust solutions that decrease possible work overloads while performing a high rate of reinsertions.

\section{Conclusion}\label{section: reinsertion-conclusion}
This paper studied MMS problem with stochastic failures and integrated reinsertion process. This is the first study that integrates the reinsertion process into a sequencing problem, to the best of our knowledge. We formulated the proposed problem as a two-stage stochastic program and proposed formulation improvements. Three bi-objective optimization algorithms were presented to tackle the problem.
The numerical experiments showed that while the two-stage bi-objective local search algorithm provides reliable solutions in terms of work overload objective, the hybrid local search integrated evolutionary optimization algorithm provides a better exploration of solution space.
To assess the quality of the solutions, dynamic reinsertion simulations are executed over industry-inspired instances. The results show that we can reduce the work overload by around 20\% while decreasing the waiting time of the failed vehicles drastically. 

\subsection{Managerial Insights}\label{section: reinsertion-managerial}
A planned sequence is disrupted first by the failure of vehicles, then by the reinsertion of the failed vehicles. In this study, we focused on generating robust schedules considering the vehicle failures and the reinsertion process which have an escalating impact on the production efficiency of the assembly lines as the difference between product types expands, currently due to electric vehicles, yet could be any other new development that induces similar results. 

The proposed problem offers the potential to generate high-quality solutions that generate a balance between minimizing work overload and waiting time of failed vehicles. This makes the problem a promising approach for real-world operational scenarios where reducing work overload and ensuring efficient reinsertions are critical considerations. By considering both failures and reinsertions, the problem addresses the inherent complexities of operations management, leading to improved performance and enhanced robustness in practical settings.

As aforementioned, the work overload reduction of around 20\% results in significant production efficiency and cost-saving enhancements in the assembly line. This is because there will be fewer line stoppages and a decrease in the number of utility workers needed. Additionally, reducing the waiting time of failed vehicles means an increase in the reinsertion ratio, which benefits decision-makers in two important ways. 
First, shorter waiting times lead to higher customer satisfaction as the late delivery of failed vehicles decreases. Second, car manufacturers impose a limit on the number of reinstating vehicles, as having a large number of piled-up vehicles disrupts the production line, such as inventory management or buffer issues. Thus, when the number of reinstating vehicles reaches the limit, the standard production stops, and only the reinstating vehicles are produced, resulting in a significant decline in production efficiency. Therefore, an increased reinsertion ratio improves production efficiency by decreasing the likelihood of such a scenario occurring.

\begin{spacing}{0.0}
\bibliographystyle{elsarticle-harv} 
\bibliography{cas-refs}
\end{spacing}

\newpage
\appendix
\section{Attainment Surface Comparison of NSGA-II, STMLS, and LS-NSGA-II} \label{appendix: reinsertion-attainment-surface}
\begin{figure}[h]
\captionsetup{justification=centering}
    \centerline{\includegraphics[width=1.0\textwidth]{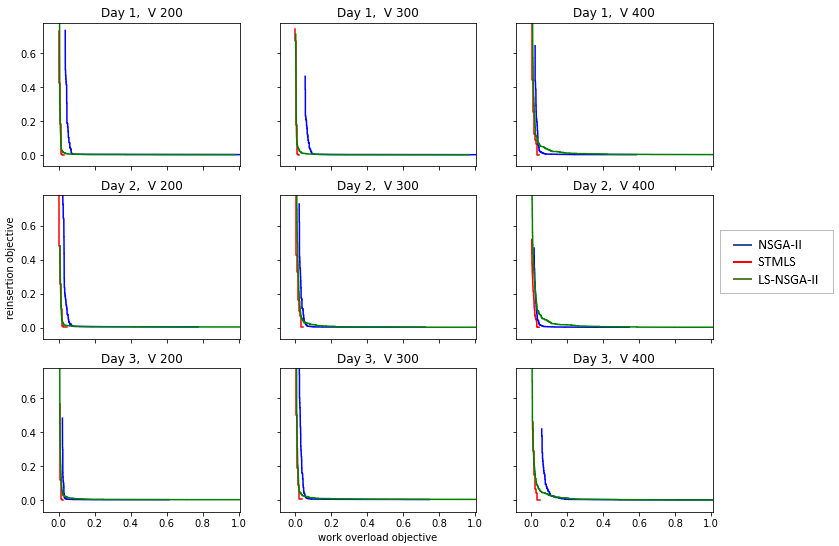}}
    \caption{Comparison of the 50\%-attainment surfaces of the proposed algorithms, with a 95\% confidence level}\label{fgr: reinsertion-attainment-global}
\end{figure}

In this section, we compare the attainment surface of the algorithms proposed in this study. Figure \ref{fgr: reinsertion-attainment-global}, which is the global counterpart of Figure \ref{fgr: reinsertion-attainment-local}, demonstrates that the LS-NSGA-II is superior to the other algorithms in terms of exploration ability. 

\end{document}